\documentclass[12pt]{article}
\makeatletter

\usepackage{amsfonts,amsmath,amssymb,theorem,graphicx}
\usepackage{amsmath,amssymb,mathrsfs,theorem,subfiles,bm}
\usepackage{graphicx}
\usepackage{subfigure}
\usepackage[top=30truemm, bottom=30truemm,left=25truemm,right=25truemm]{geometry}
\usepackage{tikz} 
\usetikzlibrary{shapes.geometric, arrows} 
\usepackage{capt-of}
\usepackage[thinlines]{easytable}
\usepackage{pbox}
\usepackage{marvosym}
\usepackage{comment}
\usepackage[colorlinks,allcolors=blue]{hyperref}  
\usepackage{here}

\newcommand{\E}{\mathbb{E}}

\def\t0{0}

\def\ia{I_{a}}
\def\is{I_{s}}

\def\p{{\mathfrak p}}

\def\rt{R_{t}}

\def\Od{\mathsf{Od}}
\def\TN{\mathsf{TN}}
\def\MN{\mathsf{MN}}
\def\tc{k}

\usepackage[normalem]{ulem}
\definecolor{c_del}{rgb}{0.8,0.1,0.1}
\definecolor{c_add}{rgb}{0.1,0.7,0.1}
\definecolor{c_cmt}{rgb}{0.1,0.1,0.7}

\title{Agent-based model and data assimilation: \\
Analysis of COVID-19 in Tokyo}

\begin{document}
\tikzstyle{block} = [rectangle, draw, thick = 1,
text width=2.5em, text centered, rounded corners, minimum height=2.5em]
\tikzstyle{block2} = [rectangle, draw, thick = 1,
text width=5.5em, text centered, rounded corners, minimum height=2.5em]
\tikzstyle{block1} = [trapezium, trapezium left angle = 65,
trapezium right angle = 115, draw, thick = 1,
text width=3em, text centered, rounded corners, minimum height=2.5em]
\tikzstyle{line} = [draw, thick = 1, -latex']
\tikzstyle{cloud} = [draw, ellipse,fill=red!20, node distance=3cm,
minimum height=2em]

\author{C. Sun${}^{a,b}$, S. Richard${}^{a,c,}\footnote{Supported by the grant \emph{Topological
                invariants through scattering theory and noncommutative geometry} from Nagoya
            University, and by JSPS Grant-in-Aid for scientific research C no 18K03328 \&
            21K03292, and on leave of absence from Univ.~Lyon, Universit\'e Claude Bernard Lyon 1,
            CNRS UMR 5208, Institut Camille Jordan, 43 blvd. du 11 novembre 1918, F-69622
            Villeurbanne cedex, France.}$\ \ \Letter, T. Miyoshi${}^{a,d,e}$}

\date{\today}
\maketitle
\vspace{-1cm}

\begin{quote}
    \emph{
        \begin{enumerate}
            \item[a)] Data Assimilation Research Team, RIKEN Center for Computational Science (R-CCS), Japan
            \item[b)] School of Science, Nagoya University, Chikusa-ku, Nagoya 464-8602, Japan
            \item[c)] Graduate School of Mathematics, Nagoya University,
                Chikusa-ku, Nagoya 464-8602, Japan
            \item[d)] Prediction Science Laboratory, RIKEN Cluster for Pioneering Research, Japan
            \item[e)] RIKEN interdisciplinary Theoretical and Mathematical Sciences Program (iTHEMS), Japan
            \item[] \emph{E-mails:} chang.sun@a.riken.jp, richard@math.nagoya-u.ac.jp, \\ takemasa.miyoshi@riken.jp
        \end{enumerate}
    }
\end{quote}

\begin{abstract}
    In this paper we introduce an agent-based model together with a particle filter
    approach for studying the spread of COVID-19. Investigations are performed
    on the metropolis of Tokyo, but other cities, regions or countries could
    have been equally chosen.
    A novel method for evaluating the effective reproduction number is one of
    the main outcome of our approach. Other unknown parameters and unknown
    populations are also evaluated.
    Uncertain quantities, as for example the ratio of symptomatic\;\!/\;\!asymptomatic
    agents, are tested and discussed, and the stability of our computations is examined.
    Detailed explanations are provided
    for the model and for the assimilation process.
\end{abstract}

\textbf{2010 Mathematics Subject Classification:}  92D30

\smallskip

\textbf{Keywords:} Agent-based model, data assimilation, COVID-19, effective reproduction number


\section{Introduction}\label{sec_intro}

The outbreak of COVID-19 had a huge impact on human society, and has
also triggered an enormous scientific response.
In this paper, we provide a rather detailed account of an agent-based model for the
evolution of the disease involving data assimilation.
Our data assimilation approach (particle filter) consists of two main steps,
repeated on a regular basis:
1) Generate a prediction for the model's state on the next time step based on the current model's state, and move the time step forward;
2) Update the model's state based on observations on the current time step.
Since new data about COVID-19 are available on a daily basis,
we regard one day as one time step in our model.
In one experiment, numerous independent simulations are considered simultaneously,
and consequently some quantities of interest can be estimated
with a probability distribution. Examples of such quantities are
unknown populations and unknown medical parameters.

The use of data assimilation for epidemiological investigations is not new, and 
several earlier references will be briefly mentioned below. 
However, up to the best of our knowledge, data assimilation has never
been paired with an agent based model in epidemiology.
For this reason and for completeness, we describe thoroughly the model and
the assimilation process. Another important aim is to make our approach and results
easily reproducible. Since some parameters in the model are location dependent,
and in order to rely on a constant source of observations,
our investigations are focusing on the metropolis of Tokyo.
Though, the approach can be adapted to other cities, regions, or countries with little effort.

The \emph{effective reproduction number} $\rt$ is one quantity of major importance
for any epidemiologic studies.
Since it works as an estimator of the number of secondary cases produced by a primary case at time $t$, it provides crucial information about the spread of the disease, allowing health services to evaluate previous interventions and forecast the evolution of the epidemic.
However, its precise evaluation is very challenging,
and several independent approaches exist
or have been recently developed. As thoroughly investigated in
\cite{GMc}, all these techniques suffer from some limitations or biases.
In this context, our approach provides one new technique for the evaluation
of $\rt$, with simple explanations about its definition and its estimation.

\begin{figure}[hbtp]
    \centering
    \includegraphics[width=12cm]{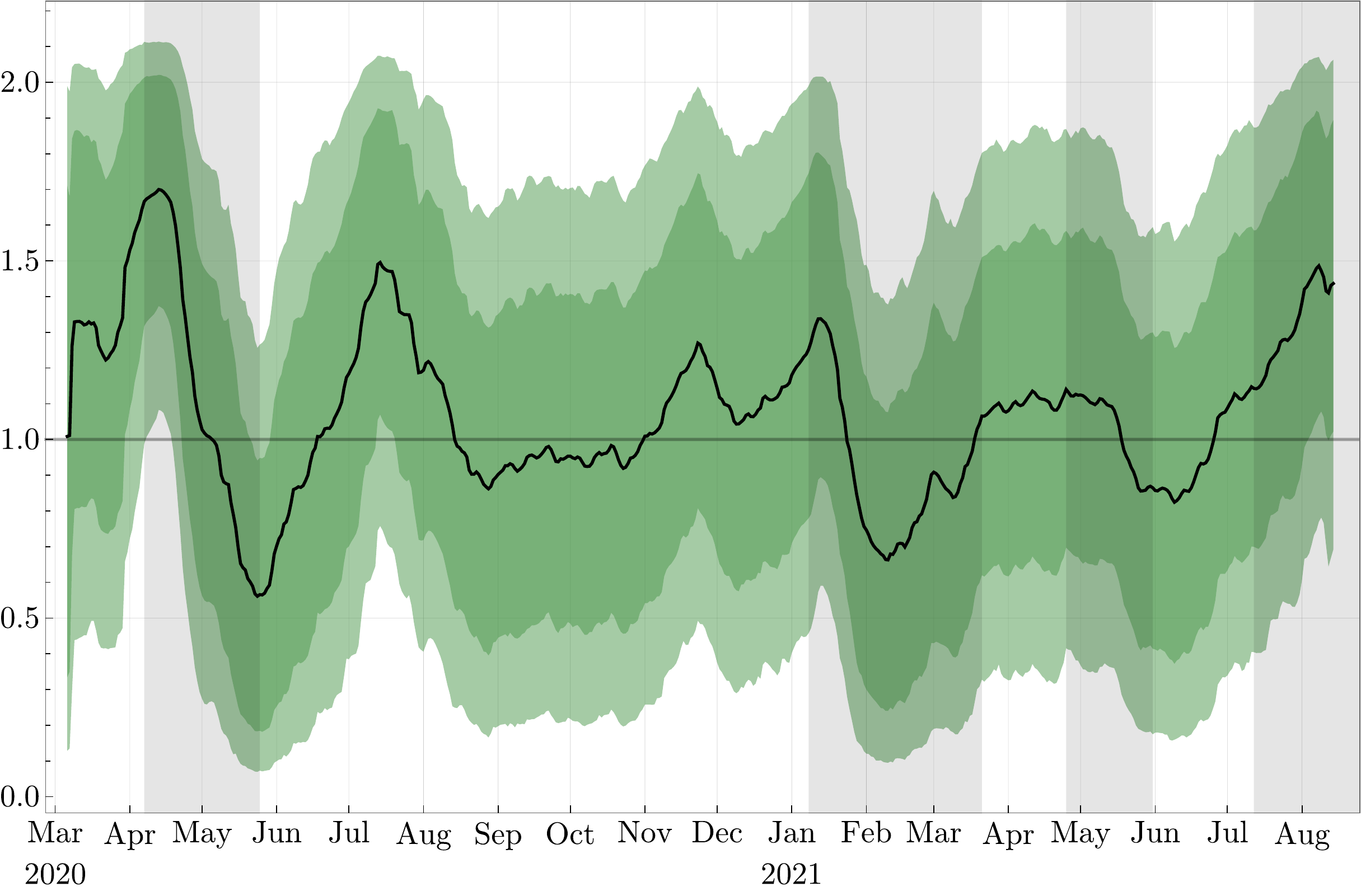}
    \caption{{\it Effective reproduction number.
                The black curve indicates the mean value; the dark and light colored regions refer to the $68\%$ and $90\%$ confidence intervals, respectively. States of emergency correspond to grey regions}}
    \label{fig_Rt}
\end{figure}

In Figure \ref{fig_Rt} we present our main result about the evaluation of the effective reproduction number $\rt$ for Tokyo, from March $6^{\rm th}$, 2020, to August $14^{\rm th}$, 2021.
Our approach does not only provide a daily mean value but
also confidence intervals based on the distribution of parameters from the individual simulations. In Section \ref{sec_dependence} we also
discuss the stability of the mean value of $\rt$
with respect to the change of several uncertain parameters.
A comparison with other estimations of the effective reproduction number for
the metropolis of Tokyo is also discussed in Section \ref{sec_final}.

As mentioned above, data assimilation has already been employed for
epidemiological investigations. One of the earliest attempts to implement epidemiological observations with epidemic models through a variational data assimilation approach was published in \cite{RH}. In relation to COVID-19, already in March 2020 some
investigations using a method called Ensemble Smoother with Multiple Data Assimilation
(ES-MDA) appeared in a preprint format \cite{vW}, quickly followed by
other publications \cite{Arm, EN, LZL, NA}. An important work in this framework
is the international collaboration \cite{E},
in which it clearly appears that the regional character of some parameters are important.
More recently, let us mention \cite{GH} in which a model with vaccination
and using an Ensemble Kalman Filter is implemented for Saudi Arabia,
and \cite{Cheng} which discusses vaccination strategies by using data assimilation techniques.
Reference \cite{AM} tracks the effective reproduction number, and
provides estimates for the effectiveness of non-pharmaceutical interventions,
while \cite{Daza} uses Bayesian sequential data assimilation for forecasting
the evolution of COVID-19 in several Mexican localities.
Additional recent investigations mixing data assimilation techniques and epidemic models
can be found in \cite{MI,RE,Sil}.

Let us now describe the content of our paper.
In Section \ref{sec_model} we precisely introduce our agent-based compartmental model (an extended version of the SEIR model)\;\!: its compartments,
the flows between these compartments, the probability of agents entering each flow, and the time agents spend in each compartment. Some of these values are fixed, in which case
we provide a reference; other values are random variables, whose
distributions are either provided by some references or will be evaluated during the experiments.
The procedure for producing secondary cases is also described, and the computation
of the effective reproduction number $\rt$ is explained.

In Section \ref{sec_evol} the evolution process and the data assimilation technique
are thoroughly introduced. 
The leading idea is to consider simultaneously numerous independent epidemics
(called particles), and to associate to each of them a weight on a daily basis.
The weights are assigned based on the proximity of the simulated epidemics to the observations. More precisely, since direct observations are only possible
for three compartments of the model (the agents under medical treatment,
the agents who have recovered after medical treatment, and the deceased agents),
we compare the number of agents in these three compartments to the corresponding
observations, and compute the weights accordingly.
Various plots illustrate the outcomes of our experiments, as for example the total
number of asymptomatic agents.
A resampling process necessary for our experiments is also introduced and discussed.
In the final part of this section, we provide some estimations about
two parameters evaluated during the experiments.

A discussion about the ratio of asymptomatic agents and the relative infectivity of symptomatic\;\!/\;\!asymptomatic agents, is presented in Section \ref{sec_dependence}. Indeed, the ratio of symptomatic\;\!/\;\!asymptomatic
agents for COVID-19 is still not precisely known, and contradictory ratios can be
found in the literature, see for example \cite{BEC, He}. 
For this reason, we have tested our model with different
ratios, and compared the resulting effective reproduction numbers.
Similarly, it is not certain by how much asymptomatic agents are less infectious than
symptomatic agents. Again, various numbers can be found in the literature \cite{BCB}, 
and for this reason we have tested our model for different relative infectivity.
At a more technical level, we also discuss in this section
the relation between the stability of the computations and the frequency of resamplings introduced in Section \ref{sec_evol}.

In the last section of this paper, we compare our results
with similar results about Tokyo, and discuss the strength and the weakness of our
approach. Some possible improvements are also presented, and future extensions
are finally discussed.

\section{Description of the model}\label{sec_model}

The model consists of an extension of the well-known SEIR model, but is adapted in the context of an agent-based model. It is based on seven compartments, as shown in Figure \ref{pesir}, which can be described as follows:
\begin{enumerate}
    \item[$S$] The compartment of all susceptible agents, they do not play any role yet.
        The size of this compartment is irrelevant for the model, and thus a few parameters related to this compartment hold an index $\infty$.
    \item[$E$] The exposed agents. Right after infections, agents moved to $E$, and they are not infectious in this compartment. They are not recorded by health authorities.
    \item[$\ia$] The infectious agents that are unaware of their infectious condition or not taking special precautions. Some of these agents are asymptomatic (will never develop any symptoms), while some are pre-symptomatic or still weakly symptomatic.
        The asymptomatic agents are removed from the model after recovery. Symptomatic agents will move to $\is$ as symptoms further develop.
They are not recorded by health authorities yet.
    \item[$\is$] The infectious agents are having clear symptoms and taking all necessary precautions. Most of these agents will receive medical support, but others will only
        quarantine themselves. The first cohort is moving to $H$ later, while the latter cohort will not be identified by the medical system and will be removed from the model. 
        They are not recorded by health authorities.
    \item[$H$] The hospitalized agents, and more generally, the agents are undergoing treatment. These agents are treated in hospitals, at home, or in other facilities. They are recorded by health authorities.
    \item[$D$] The deceased agents who were recorded by health authorities while in $H$.
    \item[$R$] The recovered agents who were recorded by health authorities while in $H$.
\end{enumerate}

\begin{center}
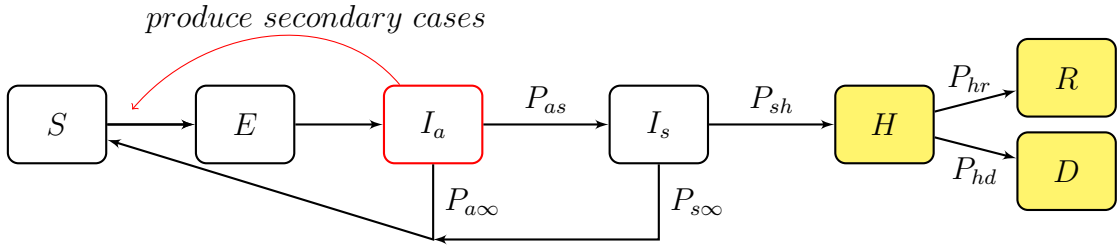

    \begin{tikzpicture}[node distance=3cm, auto]
        \node [block, draw = black] (init) {$S$};
        \node [block, draw = black, right of = init, xshift = -0.5cm] (E) {$E$};
        \node [block, draw = red, right of=E,xshift = -0.5cm] (I1) {$\ia$};
        \node [block, draw = white, above of = E, xshift = -0.8cm, yshift = -1.6cm] (note) {$produce\hbox{ }secondary\hbox{ }cases$};
        \node [block, draw = black, right of=I1] (I2) {$\is$};
        \node [block, fill=yellow!70, right of= I2] (H) {$H$};
        \node [block, fill=yellow!70, above right of= H, yshift=-1.5cm, xshift = 0.3cm] (RH) {$R$};
        \node [block, fill=yellow!70, below right of= H, yshift= 1.5cm, xshift = 0.3cm] (DH) {$D$};
        \path [line] (init) -- node {}(E);
        \path [line] (init) -- node [shift={(0,-7mm)}] {}(E);
        \path [line] (E) -- node {}(I1);
        \path [line] (I1) -- node {$P_{as}$}(I2);
        \path [line] (I2) -- node {$P_{sh}$}(H);
        \path [line] (H) -- node [near end, shift={(1mm,-1mm)}] {$P_{hr}$}(RH);
        \path [line] (H) -- node [near start, shift={(-2mm,-7mm)}] {$P_{hd}$}(DH);
        \path[line] (I1.south) |- node [near start] {$P_{a \infty}$}++(0, -10mm)--(init);
        \path[line] (I2.south) |- node [near start] {$P_{s \infty}$}++(0, -10mm)- ++(-30mm, 0);
        \draw[- >, red] (I1) to [out=130,in=50]++(-40mm, 2mm);
    \end{tikzpicture}
    \captionof{figure}{The compartments and the probabilities of each path}
    \label{pesir}
\end{center}

As shown in Figure~\ref{pesir}, should a compartment have multiple outflowing paths, its agent will enter different paths following certain probabilities. We provide in Table \ref{table_prob} the different values and the sources. Note that these values will be discussed again with additional experiments in Section \ref{sec_dependence}.
Note also that the probabilities $P_{hd}$ and $P_{hr}$ will be evaluated during the experiments, and therefore are time ($\equiv$day) dependent.
Note also that, as the probability of self-quarantining without contacting any health authorities ($P_{s\infty}$) is not available for Tokyo, we use the result of the survey \cite{OP} conducted in Osaka.

\begin{table}[htbp]
    \centering
    \begin{TAB}(r,0.5cm,0.5cm)[6pt]{|c|c|l|c|c|}{|c|c|c|c|c|c|c|}
        Name & Path & Probability of & value (\%)& source \\
        $P_{as} $ & $\ia \to \is$ & being symptomatic & 83 & \cite{BCB} \\
        $P_{a\infty}$ & $\ia \to S$ & being asymptomatic & 17 & \cite{BCB}\\
        $P_{sh}$ & $\is \to H$ & contacting health authorities & 78 & \cite{OP} \\
        $P_{s\infty}$ & $\is \to S$ & not contacting health authorities & 22 & \cite{OP}\\
        $P_{hd}$ & $H \to D$ & dying & Evaluated with observations & \\
        $P_{hr}$ & $H \to R$ & recovering after hospitalization & Evaluated with observations & \\
    \end{TAB}
    \caption{Probabilities between compartments}
    \label{table_prob}
\end{table}

In some of these compartments (Cpt), the number of days spent by agents is important.
We list in Table \ref{table_time} the necessary information about these
durations and provide the sources. Some given distributions are provided
in Table \ref{Distribution times}. Note also that
the time duration $T_h$ (time in compartment $H$) will be evaluated during the
 experiments, and therefore is time dependent.
Note that if any of the parameters are time dependent, the behavior of agents are determined by the parameters on the day they are enter their current compartments.

\begin{table}[htbp]
    \centering
    \begin{TAB}(r,0.5cm,0.5cm)[4pt]{|c|c|l|l|l|}{|c|c|c|c|c|}
        Cpt & Name & Description & Value & \hfil source \\
        $E$ & $T_e$ & Incubation period & 3 days & \pbox{2.7cm}{Estimated period before becoming contagious}\\
        $\ia$ & \pbox{6cm}{$T_{as}$ \\ $T_{a\infty}$} & \pbox{6cm}{Time before moving to $\is$ \\ Time before recovering for \\ asymptomatic} & \pbox{6cm}{Given in Table \ref{Distribution times} \\ Given in Table \ref{Distribution times}} & \cite{PL}\\
        $\is$ & \pbox{6cm}{$T_{sh}$ \\ $T_{s\infty}$} & \pbox{6cm}{Time before moving to hospital \\ Time before recovering} & \pbox{6cm}{Given in Table \ref{Distribution times} \\ Irrelevant for the model}&
        \cite{M3}\\
        $H$ &$T_h$ & Time spent in hospital & Evaluated with observations & \\
    \end{TAB}
    \caption{Time spent in compartments}
    \label{table_time}
\end{table}

\begin{table}[htbp]
    \begin{center}
        \begin{tabular}{|l|c|c|c|c|c|c|c|c|c|c|}
            \hline
            days          & 1 & 2   & 3   & 4    & 5    & 6    & 7   & 8   & 9    & $\E$ \\ \hline
            $T_{a\infty}$ & 0 & 0   & 0   & 0    & 0.05 & 0.2  & 0.5 & 0.2 & 0.05 & 7    \\ \hline
            $T_{as}$      & 0 & 0.3 & 0.6 & 0.05 & 0.05 & 0    & 0   & 0   & 0    & 2.85 \\ \hline
            $T_{sh}$      & 0 & 0.1 & 0.2 & 0.5  & 0.15 & 0.05 & 0   & 0   & 0    & 3.85 \\ \hline
        \end{tabular}
    \end{center}
    \caption{Distribution for $T_{ar}$, $T_{as}$, and $T_{sh}$}
    \label{Distribution times}
\end{table}

In compartment $\ia$, asymptotic agents and symptomatic agents can spread the disease. A simplified daily offspring distribution $\Od$ is provided in Table \ref{Offspring}, and we refer for example to \cite{Xu} for more details.

\begin{table}[htbp]
    \begin{center}
        \begin{tabular}{|c|c|c|c|c|c|c|c|}
            \hline
            $\#$ of s.~c. & 0         & 1          & 2          & 3          & 4          & 5          & $\E(\mathrm{s.~c.})$ \\ \hline
            Prob.         & $p_0=0.5$ & $p_1=0.35$ & $p_2=0.12$ & $p_3=0.01$ & $p_4=0.01$ & $p_5=0.01$ & $0.71$               \\ \hline
        \end{tabular}
    \end{center}
    \caption{Daily offspring distribution $\Od$: number of secondary cases (s.c.)}
    \label{Offspring}
\end{table}

One delicate question is the relation between the transmission coefficient for asymptomatic agents and the transmission coefficient for symptomatic agents. For our investigations, we shall rely on the result of the systematic review \cite{BCB} which asserts that the relative risk of asymptomatic transmission
is $42\%$ lower than that for symptomatic transmission. As a consequence, we shall fix
that the \emph{relative infectivity coefficient $\tc$} is $0.58$, which scales the asymptomatic agents' transmission coefficient relative to the symptomatic agents'.
This factor is slightly smaller but of a comparable scale compared to earlier investigations, see for example \cite{BEC}. This factor will be discussed again with additional experiments in Section \ref{sec_dependence}.

Let us set $r_t\in [0,1]$ for a time dependent multiplicative factor which takes into account the real interaction between agents.
This factor (updated on a daily basis) depends clearly on individual behavior but also on non-pharmaceutical interventions.
As a consequence, for asymptomatic agents in $\ia$, the daily production of second generation infections is given by the daily offspring production with all $\#\;\!$secondary cases $\geq 1$ entries' probabilities multiplied by the factor $\tc \cdot r_t$. For
symptomatic agents in $\ia$, the daily production of second generation infections is given by the same rule, but only with multiplicative factor $r_t$ instead of $\tc \cdot r_t$.
Therefore, the effective reproduction number $R_t$ is given by the formula
\begin{align}\label{eq_rep}
    R_t       & = \left(\tc\cdot P_{a\infty}\cdot \E\big(T_{a\infty}\big) + P_{as}\cdot\E\big(T_{as}\big) \right)\cdot \E(\mathrm{s.~c.})\cdot r_t \\
    \nonumber & = (0.58\cdot 0.17\cdot 7 + 0.83\cdot 2.85) \cdot 0.71\cdot r_t                                                                     \\
    \nonumber & = 2.17 \cdot r_t.
\end{align}

\section{Evolution process, particle filter, and main results}\label{sec_evol}

The observation data for compartments $H$, $R$, and $D$ provided by health authorities for Tokyo start on March $6^{\rm th}, 2020$. However, the epidemic had already started
in January at the latest, since the departure
of the Diamond Princess from Yokohama took place on January $20^{\rm th}, 2020$.
For this reason, and after several trials, we have fixed the following initial conditions
for our experiment: January $17^{\rm th}, 2020$\footnote{Long after having chosen this date as the most suitable one for the start of our simulations, we were informed that the first case of COVID-19 detected in Japan was on January $16^{\rm th}$, 2020.} . 
This date corresponds to $49$ days before the start of the observations available for Tokyo.
\begin{table}[htbp]
    \begin{center}
        \begin{tabular}{|l|l|}
            \hline
            Description                       & Initial values                                               \\ \hline
            Initial date                      & January 17, 2020                                             \\ \hline
            Initial number of infected agents & Uniformly at random in $\{3,4,5,6,7\}$                       \\ \hline
            $r_0$                             & Uniformly at random in $[0,1]$                               \\ \hline
            $T_h(0)$                          & $15$ days, initial guess motivated by \cite{TMG}\\ \hline
            $P_{hd}(0)$                       & Uniformly at random in $[0, 0.05]$                           \\ \hline
        \end{tabular}
    \end{center}
    \caption{Initial conditions}
    \label{Table_initial}
\end{table}

On a daily basis, some parameters have certain freedom to change their values.
In order to describe their evolutions, let us use the notation $\TN(\mu,\sigma;a,b)$
for the truncated normal distribution of mean $\mu$, variance $\sigma^2$, minimum value $a$ and maximal value $b$.
For the evolution of $r_t$ with $t\geq 1$, we choose randomly its value according to the distribution $\TN\big(r_{t-1}, 0.05; 0, 1\big)$.
The evolution of $P_{hd}$ is also given by using a truncated normal distribution,
namely, for $t\geq 1$ the value $P_{hd}(t)$ is chosen randomly according to the distribution $\TN\big(P_{hd}(t-1),0.0025; 0, 0.05\big)$

Since the time spent in hospitals, $T_h$, is integer valued, its evolution is slightly more complicated, but nevertheless follows a scheme similar to the previous two parameters.
The mean value $\E\big(T_h(t)\big)$ of $T_h(t)$ is chosen randomly according to the distribution $\TN\big(\E(T_h(t-1)), 0.75; 4,19\big)$, where the minimum and maximum values have been fixed according to some information provided by health authorities \cite{M3}. Then, one constructs a distribution supported on the greatest integer less than or equal to $\E\big(T_h(t)\big)$ and the least integer greater than or equal to $\E\big(T_h(t)\big)$, and such that its expectation is equal to $\E\big(T_h(t)\big)$. The agents entering the compartment $H$ on that day are then assigned a time in $H$ chosen at random according to this distribution.

Each agent in the compartment $\ia$ will infect susceptible agents belonging to the compartment $S$ on a daily basis. We now describe this process.
Let us denote by
$$
    \MN(x_j,n;\p_j)\equiv \MN(0,1,2,3,4,5,n; \p_0,\p_1,\p_2,\p_3,\p_4,\p_5),
$$
the multinomial distribution for $n$ trials in the set of values $x_j\in \{0,1,2,3,4,5\}$
with the probability of $x_j$ given by $\p_j$.
We also denote by the vector $X=(X_0,X_1,X_2,X_3,X_4,X_5)$ one realization of this distribution. For example, if $\p_j$ is the the probability of having $j$ offsprings as given in Table \ref{Offspring}, then one has $\sum_{j=0}^5 X_j = n$, and the expectation value for $X_j$ is $n p_j$ for any $j\in \{0,1,2,3,4,5\}$.

On a given day $t$, assume that the compartment $\ia$ contains $n$ asymptomatic agents and $m$ symptomatic agents. Then, the $n$ agents will infect $\sum_{j=1}^n j X_j$
susceptible agents, where $X$ is one realization of the multinomial distribution
$\MN(x_j,n;\p^a_j)$, with $\p^a_j = \tc \cdot r_t\cdot p_j$ for $j\in \{1,2,3,4,5\}$ and $\p^a_0 = \big(1-\tc \cdot r_t\big)+\tc \cdot r_t\cdot p_0$.
Similarly, the $m$ agents will infect $\sum_{j=1}^n j Y_j$ susceptible agents, where $Y$
is one realization of the multinomial distribution
$\MN(x_j,m;\p^s_j)$, with $\p^s_j = r_t\cdot p_j$ for $j\in \{1,2,3,4,5\}$ and $\p^s_0 = \big(1- r_t\big)+r_t\cdot p_0$.

As highlighted in Figure \ref{pesir}, three compartments $H$, $D$, and $R$, are associated with observations (collected from \cite{Toyo}). Some uncertainties must be attached to these observations, and these uncertainties are commonly called observation error and will be denoted by $\sigma$. For our experiments, we consider two of these observation errors constant over time, while one will be time dependent. In fact, since the daily variations of $H$ are quite important, the corresponding uncertainties will take them into account. More precisely, if we write $H(t)$ for the number of hospitalized agents at time $t$, we set
\begin{align*}
    \sigma_H(t) & = \sqrt{\big(0.3 H(t)+4 |H(t)-H(t-1)|\big)^2+400}, \\
    \sigma_R    & = 2000,                                            \\
    \sigma_D    & = 100.
\end{align*}
These expressions for the errors have been chosen after several trials. As shown later, even if these values look very large, they lead to a successful selection process.

Let us now describe more precisely the experimental setup.
We use a particle filter approach; namely we create a large number $N$
of independent simulations (=particles) of the propagation of the epidemic in Tokyo in one experiment.
Each of the simulations is a realization of the model, and involves all the
daily random processes described above. Typically, we have chosen $N=50'000$ or $N=100'000$, 
and kept this number constant during each experiment.
On a daily basis, a weight $w$ is assigned to each particle.
Let $i$ denote the index of the particle, and let $H_i(t)$, $R_i(t)$, and $D_i(t)$ denote the number of agents in the three corresponding compartments for this particle at time $t$.
At each time $t\geq 1$ one first computes the
not normalized weight $W_i(t)$ by the formula
$$
    W_i(t):=w_i(t-1)\cdot \exp\bigg(-\frac{(H_i(t)-H(t))^2}{2 \sigma_H(t)^2}-\frac{(R_i(t)-R(t))^2}{2 \sigma_R^2}-\frac{(D_i(t)-D(t))^2}{2 \sigma_D^2}\bigg),
$$
with the convention that $w_i(0)=\frac{1}{N}$. Then, the normalized weight for the particle $i$ is given by
$$
    w_i(t):=\frac{W_i(t)}{\sum_{j=1}^N W_j(t)}.
$$
With these weights, one gets three distributions for the analyzed values of $H$, $R$,
and $D$ at time $t$, respectively.

Now, if we keep computing weights with the formula described above,
it will soon turn out that very few particles will concentrate almost all weights, and that nearly all other particles would end up with negligible weights.
For dealing with this problem, a process called resampling is implemented:
For a given time $t$, we take one realization $X$ of the multinomial distribution $\MN\big(i,N;w_i(t-1)\big)$ with $i \in \{1,\dots,N\}$. Then, for a particle indexed $i$, we create $X_i$ copies of the particle and remove the original one. Finally, we assign a weight $1/N$ to all new particles. Clearly, some particles will appear several times, but their trajectories will diverge due to the randomness involved on a daily basis.

\begin{figure}[htbp]
    \centering
    \subfigure[]{\includegraphics[width=8.3cm]{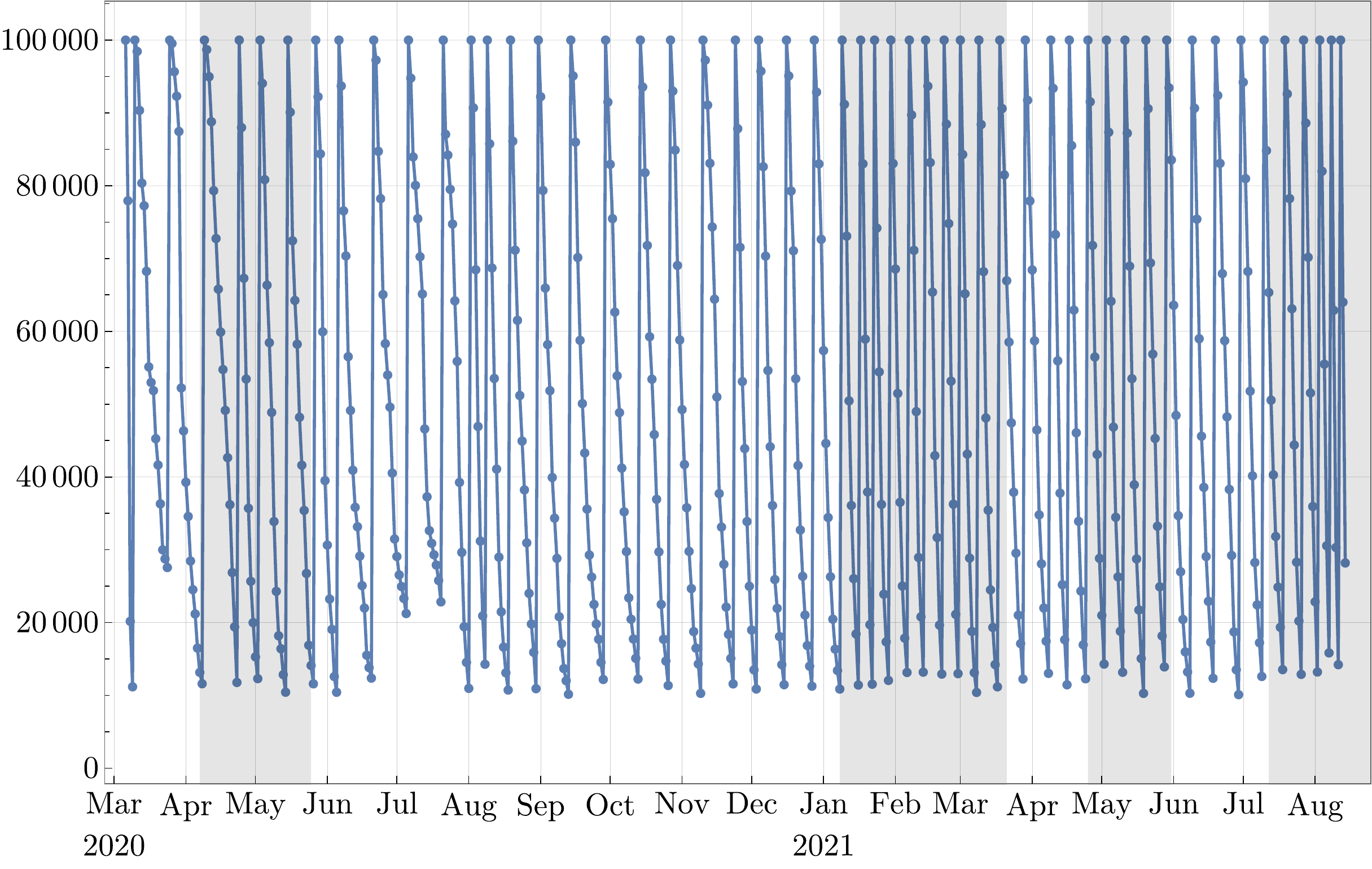}}
    \subfigure[]{\includegraphics[width=8.1cm]{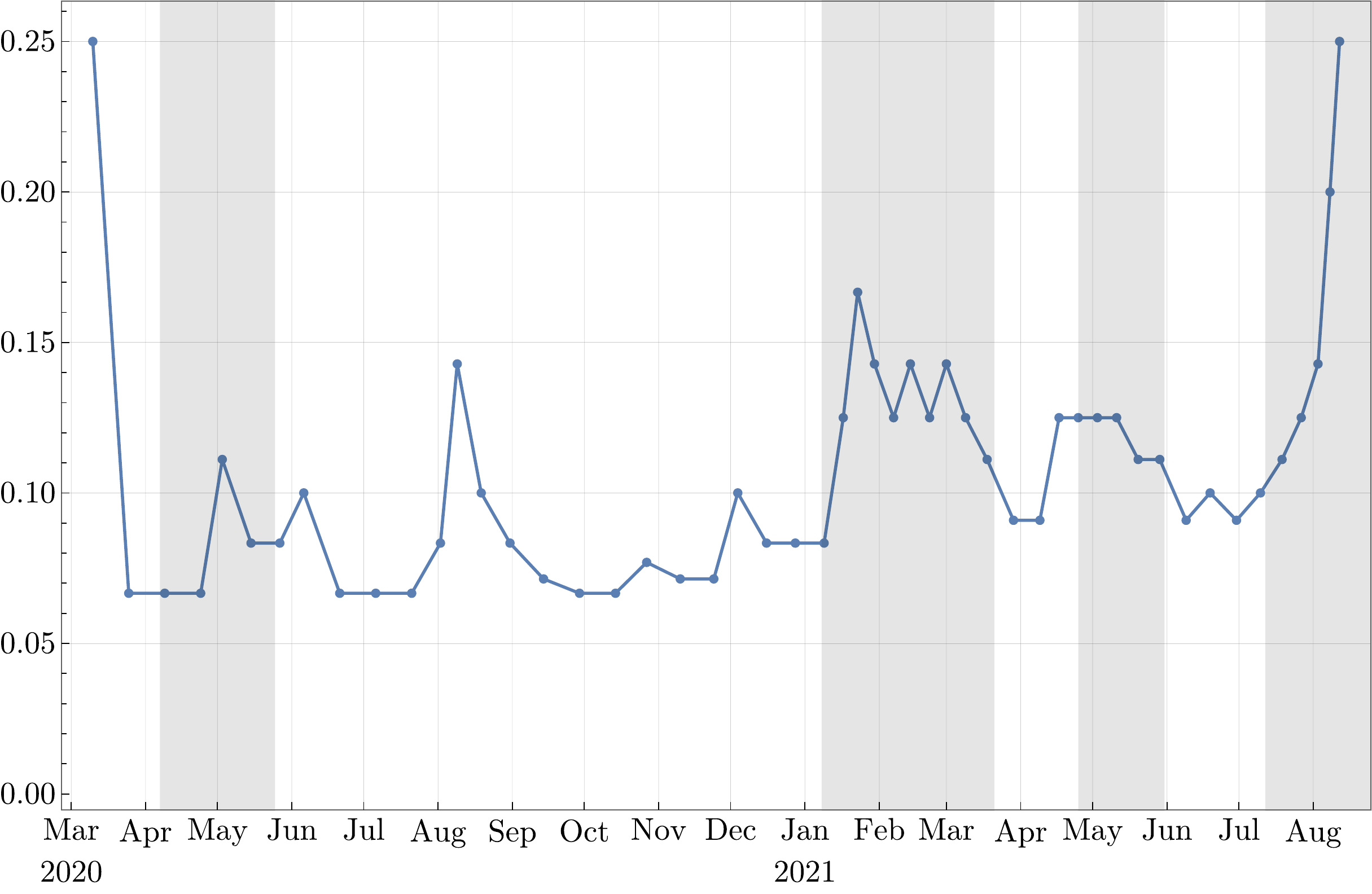}}
    \caption{{\it (a) $N_{eff}$, (b) Inverse of time between resamplings }}
    \label{fig_resampling}
\end{figure}

One still has to decide when resamplings are organized. For that purpose, let us set the effective number of particles
$$
    N_{eff}(t):=\frac{1}{\sum_{i=1}^N w_i(t)^2},
$$
and observe that if $w_i = \frac{1}{N}$ for $i\in \{1,\dots,N\}$, then $N_{eff}=N$, while if $w_i\approx 1$ for one $i$, and $w_j\approx 0$ for all $j\neq i$, then
$N_{eff} \approx 1$. For this reason, $N_{eff}$ is often used as an indicator of the number of particles still playing a role in the process. As a consequence, we shall use the following rule:
If $N_{eff}<\frac{N}{10}$, then a resampling has to take place. 
As a practice to stabilize the experiments, a resampling is also performed if the most recent resampling took place 15 days ago.
In Figure \ref{fig_resampling}(a) we present a typical graph for $N_{eff}$ for a number of particles $N=100'000$, while in Figure
\ref{fig_resampling}(b) we provide an indication about the time between the resamplings, namely whenever a resampling is taking place, we set
$y=1/(\hbox{number of days since the last resampling})$. We shall soon see that this information is playing a role for the stability of the predictions.

Let us stress that the computation of the weights and the implementation of the resampling start only on day $50$ of the simulation, namely on March 6. Indeed, as there is no observation available before this date, the particles just evolve freely and independently.

\begin{figure}[htbp]
    \centering
    \includegraphics[width=14cm]{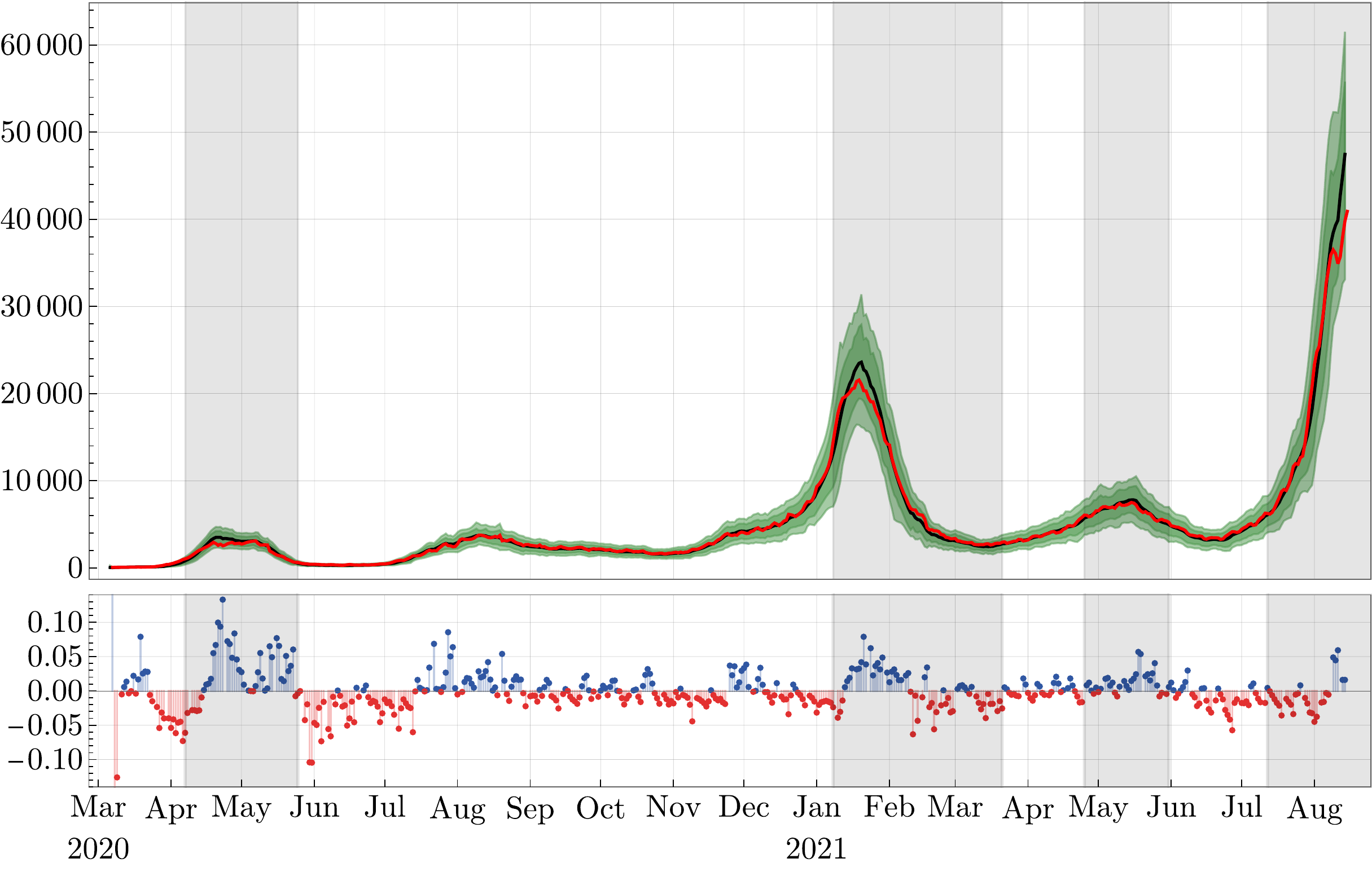}
    \caption{{\it Upper part: Observation value of $H$ (in red) and analyzed values with mean value (in black) and $68\%$, resp.~$90\%$, confidence interval. Below part: relative difference for the mean values of $1$-day forecast and analyzed $H$.}}
    \label{fig_H}
\end{figure}

After running the experiment with observation data up to day $t$, we run the experiment for one more day without observation data. With the weight computed on day $t$ still assigned to each particle, one obtains the so-called forecast (or prior) values for $H_i$, $R_i$, and $D_i$,
on day $t+1$. These values generate the corresponding forecast distributions
for $H$, $R$, and $D$. Then, by computing the weights
$w_i(t+1)$, as explained above, one obtains the analyzed values for the compartments, also called posterior values, and the corresponding analyzed distributions.
Clearly, the forecast distributions and the analyzed distributions do not match in general.
The relative difference between the mean values of the forecast and the analyzed
distributions for $H$ can be visualized in Figure \ref{fig_H} (below part).
Let us just emphasize that the weights are computed with three compartments,
not only with $H$.

\begin{figure}[h!]
    \centering
    \subfigure[]{\includegraphics[width=8.3cm]{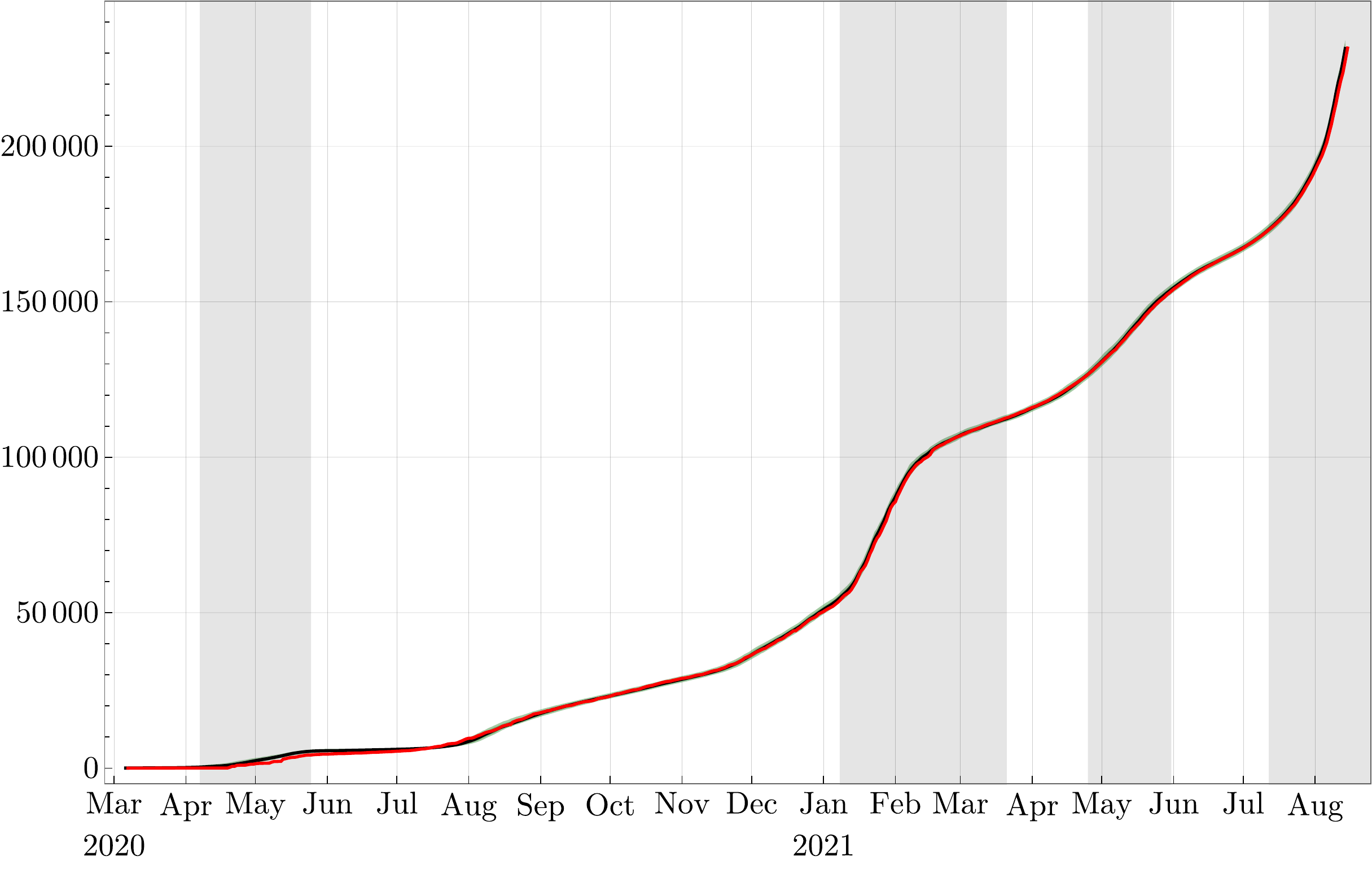}}
    \subfigure[]{\includegraphics[width=8.1cm]{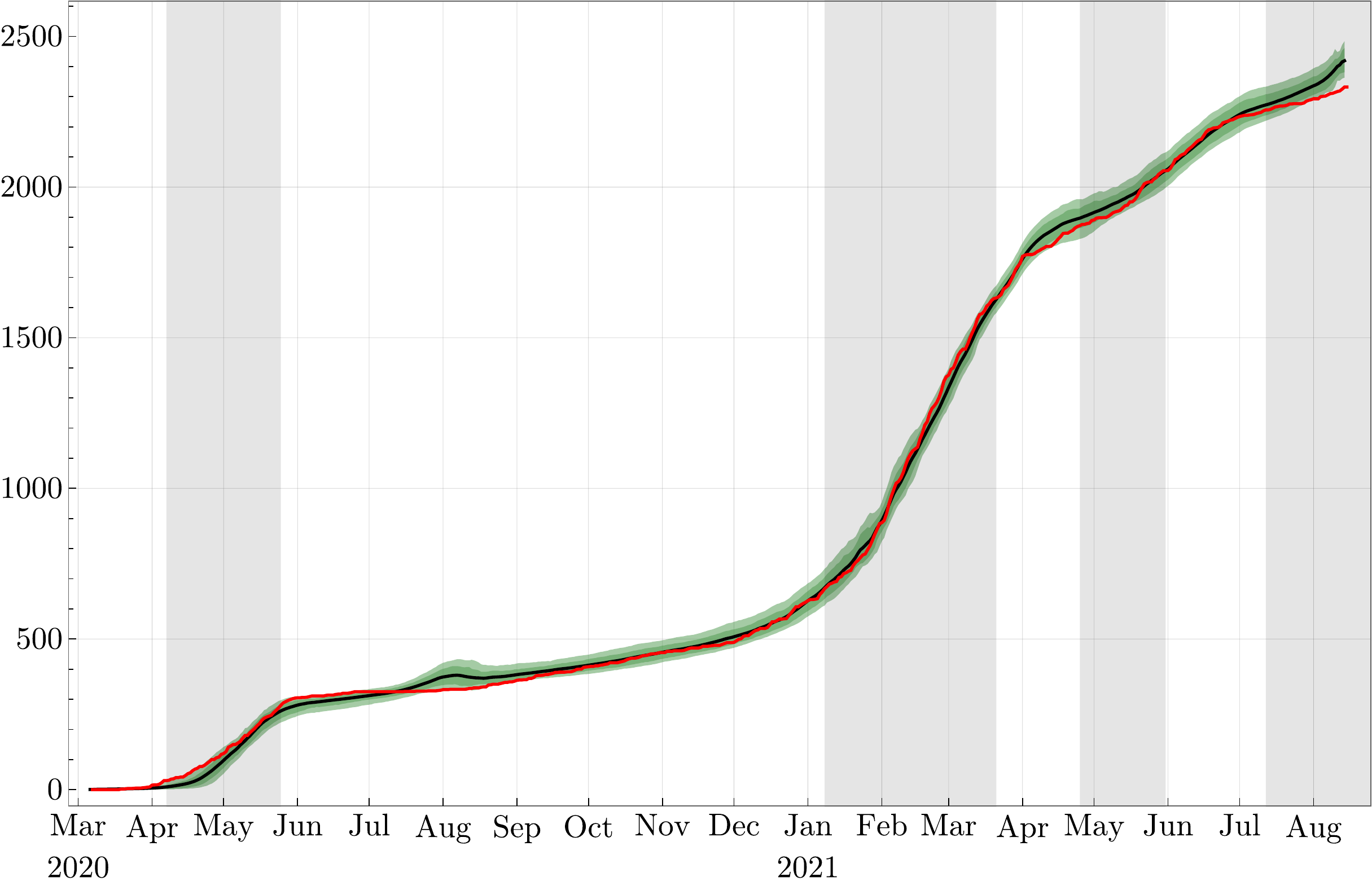}}
    \caption{{\it Observation value (in red) and analyzed values with mean value (in black) and $68\%$, resp.~$90\%$, confidence interval: (a)~Compartment $R$, (b)~Compartment $D$. In (a), all curves are superposed}}
    \label{fig_assim}
\end{figure}

\begin{figure}[htbp]
    \centering
    \includegraphics[width=8.2cm]{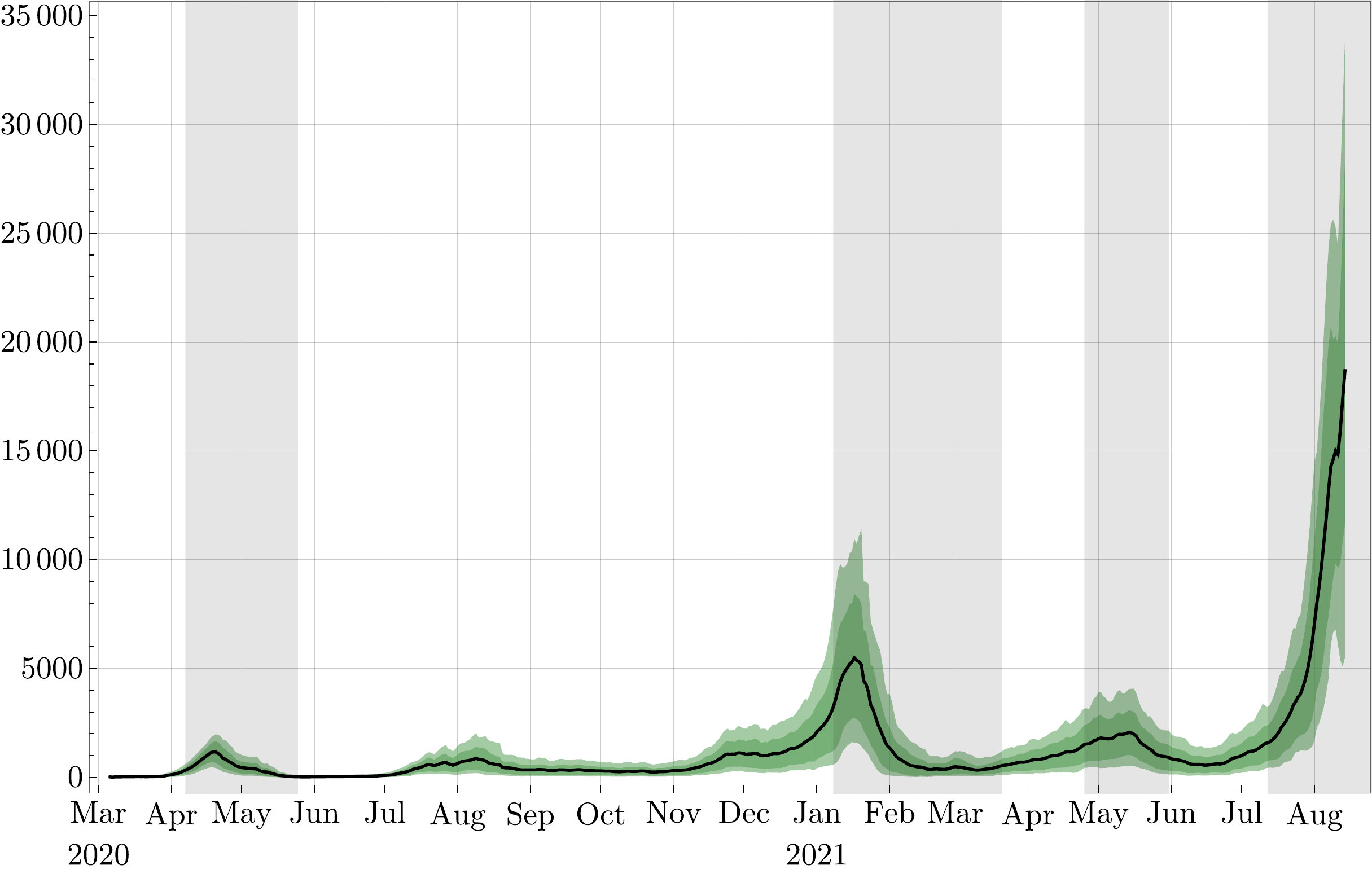}
    \caption{{\it  Total number of asymptomatic agents}}
    \label{fig_A}
\end{figure}

A representation of the three compartments with observations is also presented in Figure \ref{fig_H} (upper part) and in Figure \ref{fig_assim}.
In these plots, the observations are represented in red, and the analyzed distributions are represented in green (confidence intervals) and in black (means).
Similarly, we provide in Figure \ref{fig_A} the plot of the total number of asymptomatic agents in $\ia$.
By the definition of asymptomatic agents, there is no observation corroborating these values.

\begin{figure}[htbp]
    \centering
    \subfigure[]{\includegraphics[width=8.2cm]{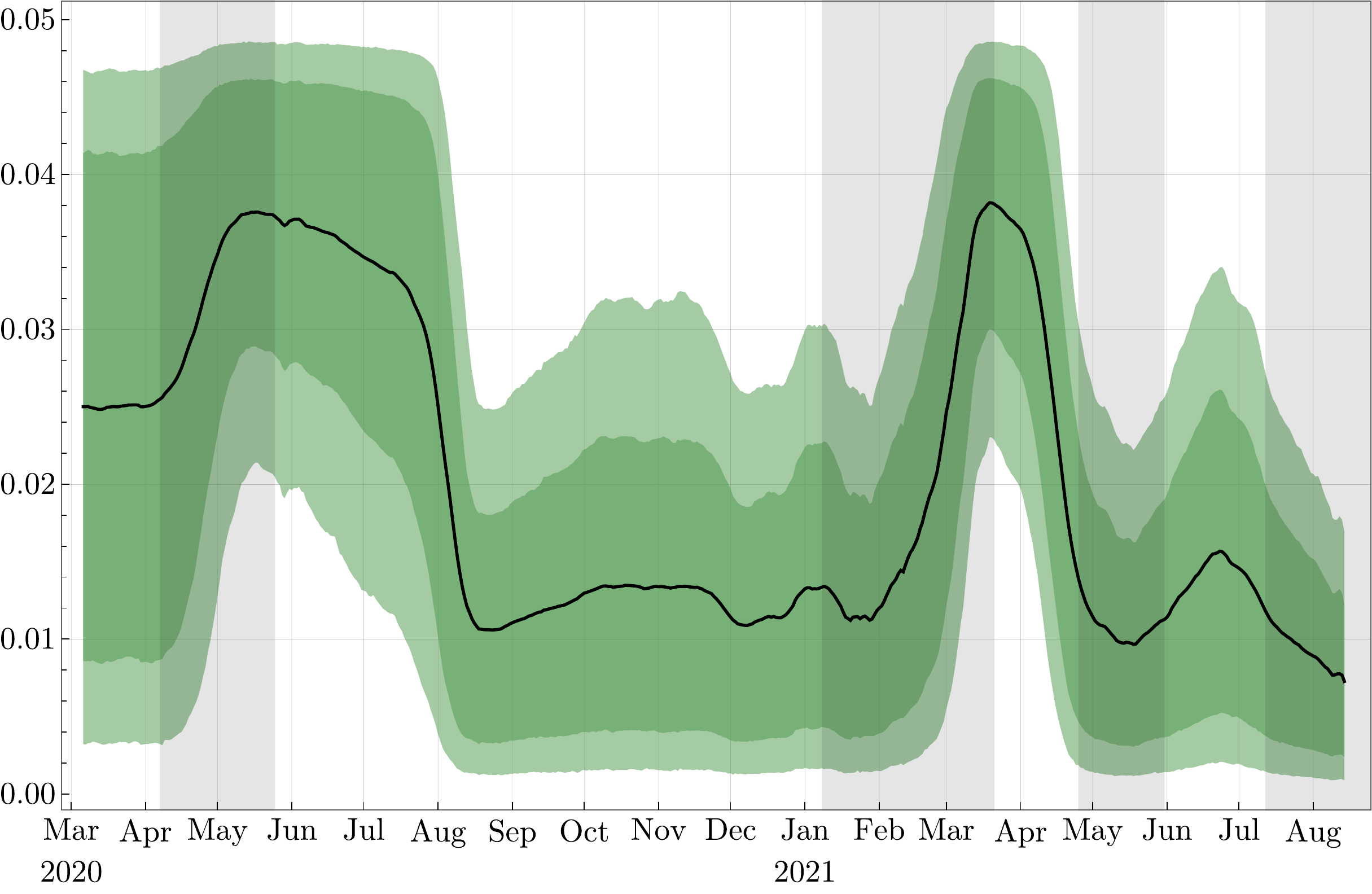}}
    \subfigure[]{\includegraphics[width=8.1cm]{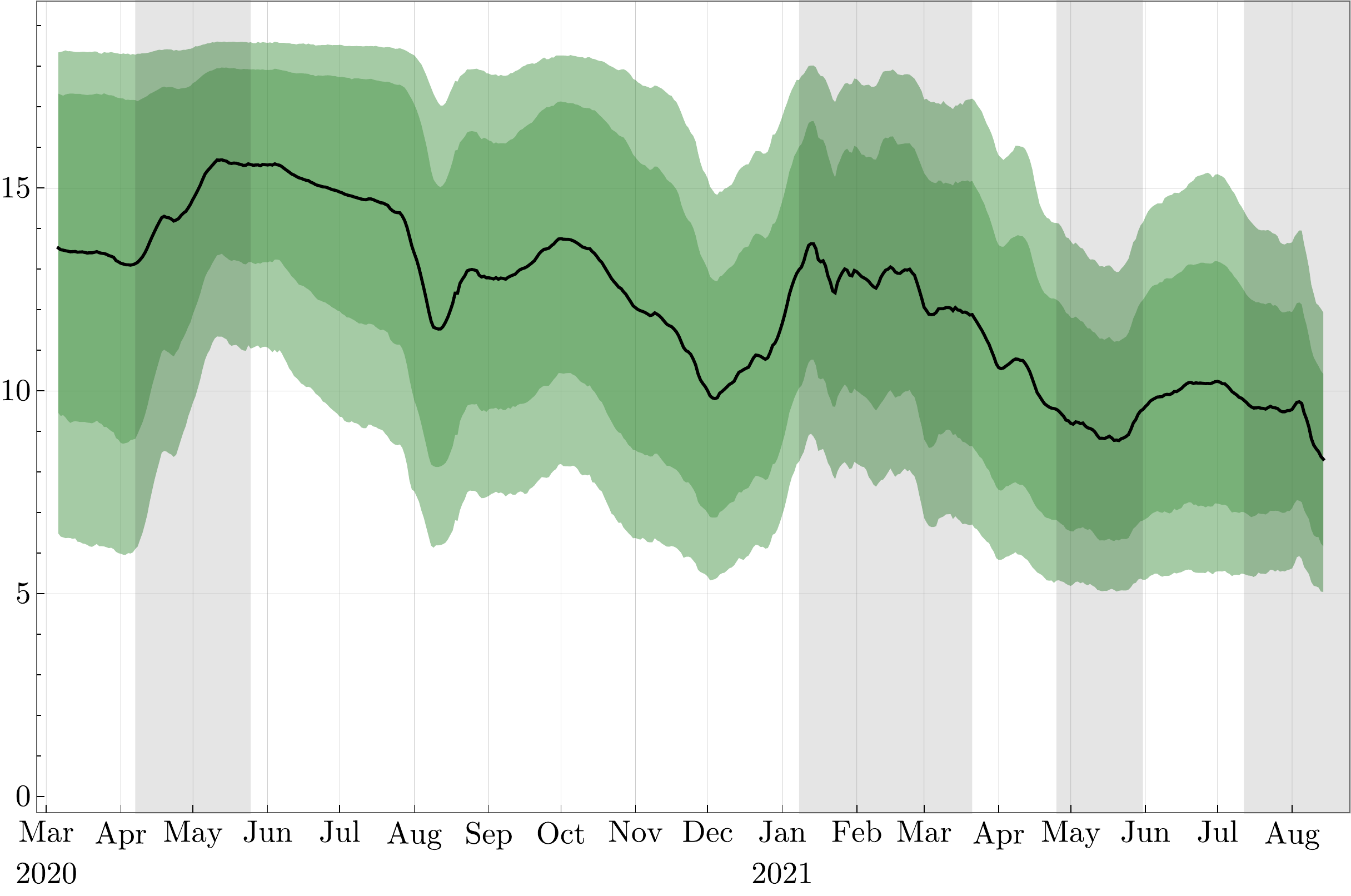}}
    \caption{{\it Mean and distribution with $68\%$ and $90\%$ confidence intervals for: (a)~$P_{hd}$, (b)~$\E(T_h)$}}
    \label{fig_death}
\end{figure}

As mentioned in Section \ref{sec_model}, $P_{hd}$, $P_{hr}$, and $\E(T_h)$ are evaluated during the experiments, which are respectively the probability of dying, the probability of recovering from $H$, and the average time agents spend in $H$.
With our approach, these quantities are provided with their distributions created by their values taken by the $N$ particles.
Clearly, if more accurate medical information were available, these parameters could be directly implemented in our model, but in their absence, we have to evaluate 
them from the observations.
For the expected value of $P_{hd}$ shown in Figure \ref{fig_death}(a), observe that the sudden increase visible between mid-February and mid-March 2021 is due to a conjunction of two factors: a large number of death, see Figure \ref{fig_assim}(b),
and the rate for the number of agents recovering which takes a local minimum, see Figure \ref{fig_assim}(a).
For the expected time $\E(T_h)$ spent in $H$, we provide in Figure
\ref{fig_death}(b) the distributions $\E\big(T_h\big)$ deduced from our investigations.
The observed long-term trend of decay is expected to be due to the improvement of the medication
and treatment for infected patients. Unfortunately, we can not deduce the difference of time spent in $H$ between patients with light symptoms and patients with serious symptoms from our model.

\section{Parameters' dependence and stability}\label{sec_dependence}

\begin{figure}[htbp]
    \centering
    \subfigure[]{\includegraphics[width=8.2cm]{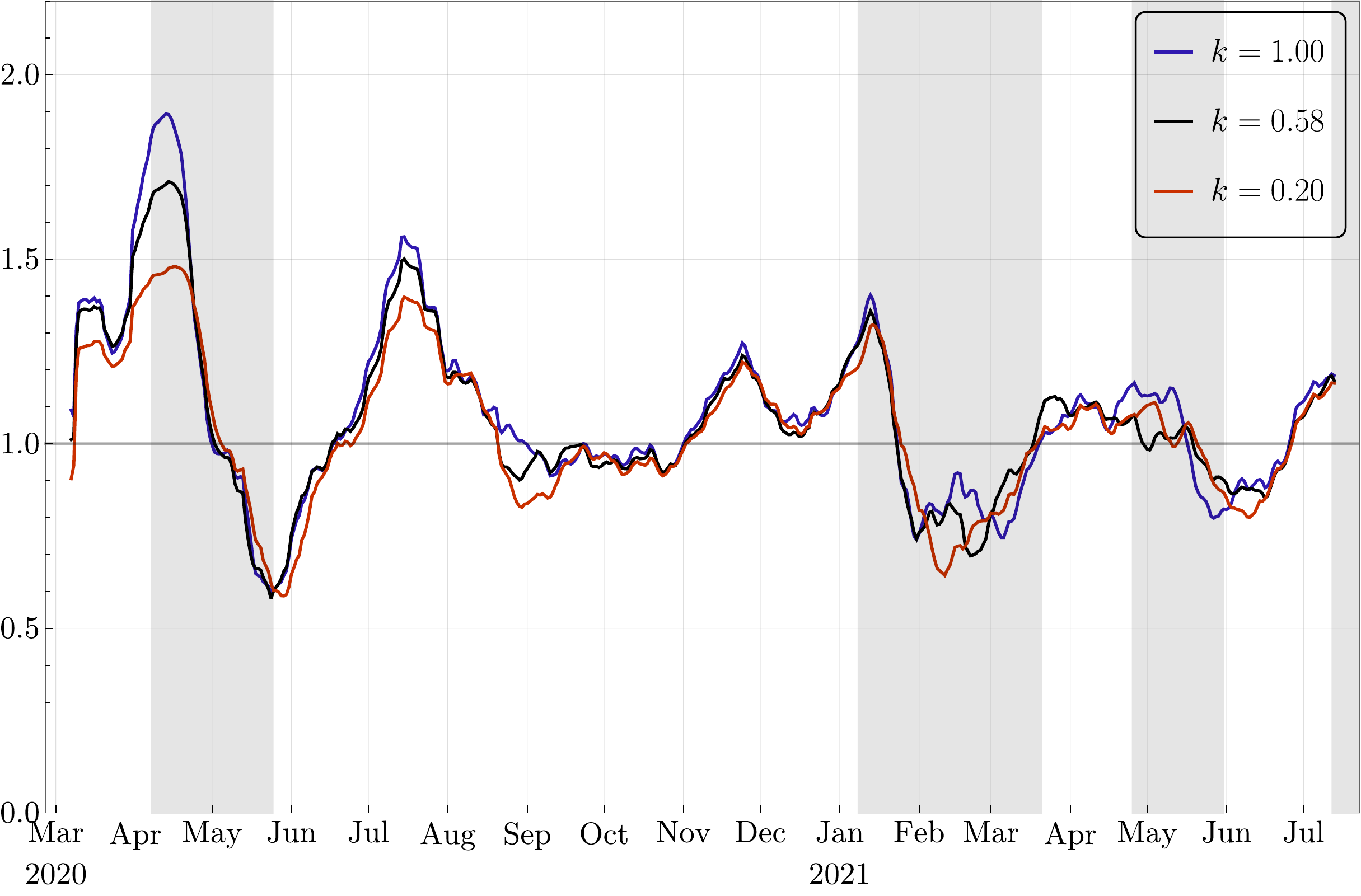}}
    \subfigure[]{\includegraphics[width=8.2cm]{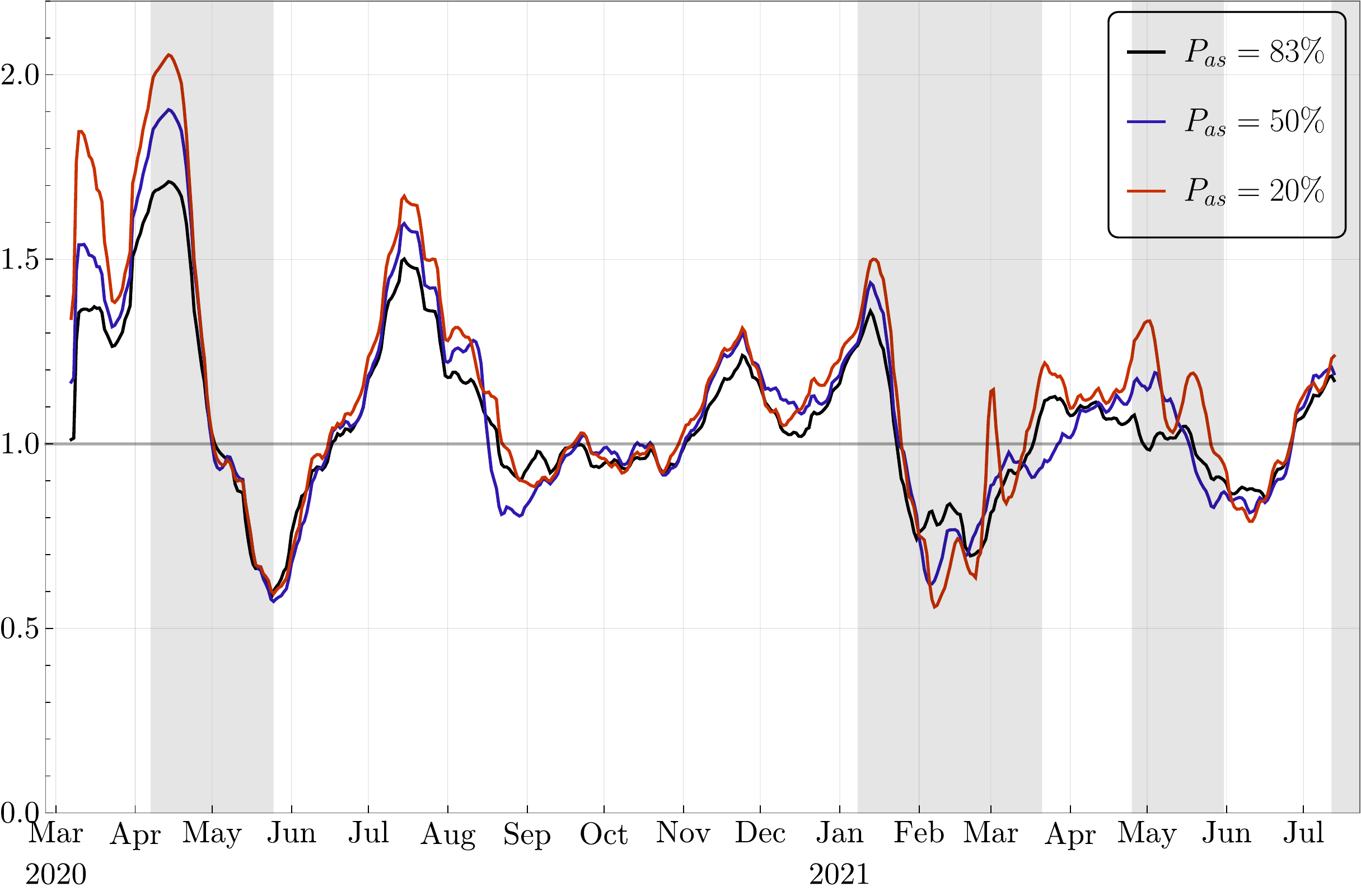}}
    \caption{{\it Mean $\rt$ for: (a)~Different values of $k$, (b)~Different ratios between symptomatic and asymptomatic}}
    \label{fig_dependence}
\end{figure}

As mentioned in Section \ref{sec_model},
the relative infectivity $k$ of the asymptomatic agents compared to the symptomatic agents
is one very uncertain parameter. For most of our simulations, we have used $k=0.58$
based on the information provided by the literature \cite{BCB}.
Since this ratio is quite uncertain, we compared the outcomes of simulations with all conditions the same but with $k=0.2$, $0.58$, and $1.0$.
The corresponding mean values for the effective reproduction number $\rt$
are shown in Figure \ref{fig_dependence}(a).
The patterns are similar, but a larger $k$ corresponds
to a slightly larger value of $\rt$, most of the time. This is not surprising since the factor $k$
plays a role in the computation of the effective reproduction number, as shown in
\eqref{eq_rep}.
On the other hand, with one noticeable exception around late August to early September 2020, the three curves cross the critical value $\rt=1$ more or less simultaneously.

Similarly, the ratio between symptomatic and asymptomatic agents is also a somewhat controversial parameter.
For our simulations, we have used the respective proportion of $83\%$ and $17\%$ from \cite{BCB}, but other sources mentioned a very different ratio \cite{He}.
Since asymptomatic cases are very difficult to detect, and since we can not be fully confident in this ratio, a sensitivity test has been performed.
To do this, we have decreased the ratio of symptomatic to $50\%$ and to $20\%$ and performed the whole experiment while keeping all other parameters the same as the original ones.
Compared to the original setting, more agents become asymptomatic and recover without showing any symptoms in these two new scenarios. On the other hand,
the number of symptomatic agents is more
or less constant since they are dominated by the number of agents hospitalized, which is compared and adjusted according to the observations on a daily basis.
In Figure \ref{fig_dependence}(b), the different curves for the mean $\rt$ have similar patterns.
However, we observe that a bigger proportion of asymptomatic agents leads to a slightly bigger $\rt$.
Indeed,  more total infected agents have to be created in this scheme, which means that a larger $\rt$ is necessary.

\begin{figure}[htbp]
    \centering
    \includegraphics[width=12cm]{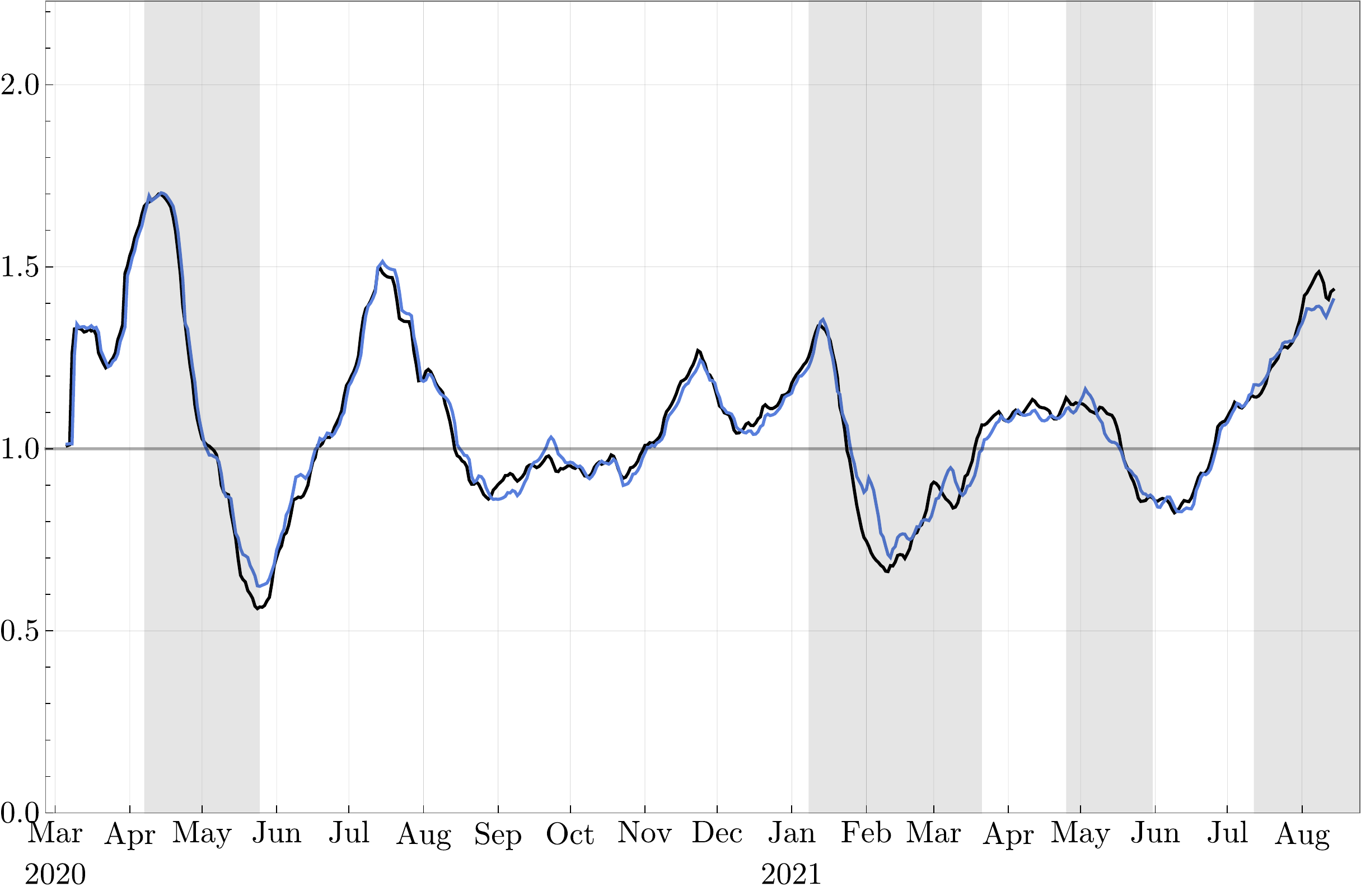}
    \caption{{\it The mean values of the effective reproduction number obtained with two independent experiments, each involving 100'000 particles}}
    \label{fig_dif}
\end{figure}

Let us still mention one observation linking the stability of the computation of $\rt$
with the frequency of resamplings.
In Figure \ref{fig_dif} we provide the mean value of $\rt$ for two independent
experiments under exactly the same setting, each of them involves 100'000 particles. Clearly, these two curves
are very similar, but at a few places, small differences are visible, as for example
in September 2020, in February and first half of March 2021, and in August 2021.
Note that the same kind of discrepancies can also be observed in Figure
\ref{fig_dependence} (around the first 2 mentioned periods).
Now, if we look back to Figure \ref{fig_resampling}(b),
it turns out that these three periods follow or coincide with periods of time
when resamplings took place more frequently (local maxima of the curve presented in Figure \ref{fig_resampling}(b)).
Our understanding is the following: when abrupt changes in the observations are taking place, and\;\!/\;\!or when the model is not accurate enough, more resamplings are necessary to keep the experiment following the observation data since only a small portion of particles will have the correct behavior. As the weights concentrate on those particles, $N_{eff}$ drops rapidly, and more resamplings are triggered as a consequence. This process also quickly eliminates the short term diversity of particles: since most of the particles are sampled from that small portion of particles, they will share similar behavior for a few days. This could lead to the following two consequences in the short term: 1) the system is less capable of adapting to further rapid changes, 2) the randomness plays a more important role.

\section{Discussion and conclusion}\label{sec_final}

The effective reproduction number $\rt$ is certainly one of the most important parameters
for the study of an epidemic, but it suffers intrinsically from several weaknesses.
Indeed, because of the delay between the infection of an agent and the appearance of symptoms, the effect of the current $\rt$ will only be visible in a few days.
Accordingly, the current situation (number of agents in $\ia$, $\is$, or in $H$)
is related to some effective reproduction numbers which took place over the last few days.
Also, since the infectious period lasts for several days, it is not possible (as for example mentioned in \cite{GMc}) to shift $\rt$ by some days to get a better picture.
These drawbacks are well-known, and are taken into account by medical institutions.
Our approach does not solve this problem, but it provides an alternative way of estimating $\rt$, and makes it clear that $\rt$
is really the reproduction number taking place on day $t$. As emphasized in \cite{GMc},
this synchronicity is not shared by all methods estimating $\rt$.

\begin{figure}[t]
    \centering
    \includegraphics[width=15cm]{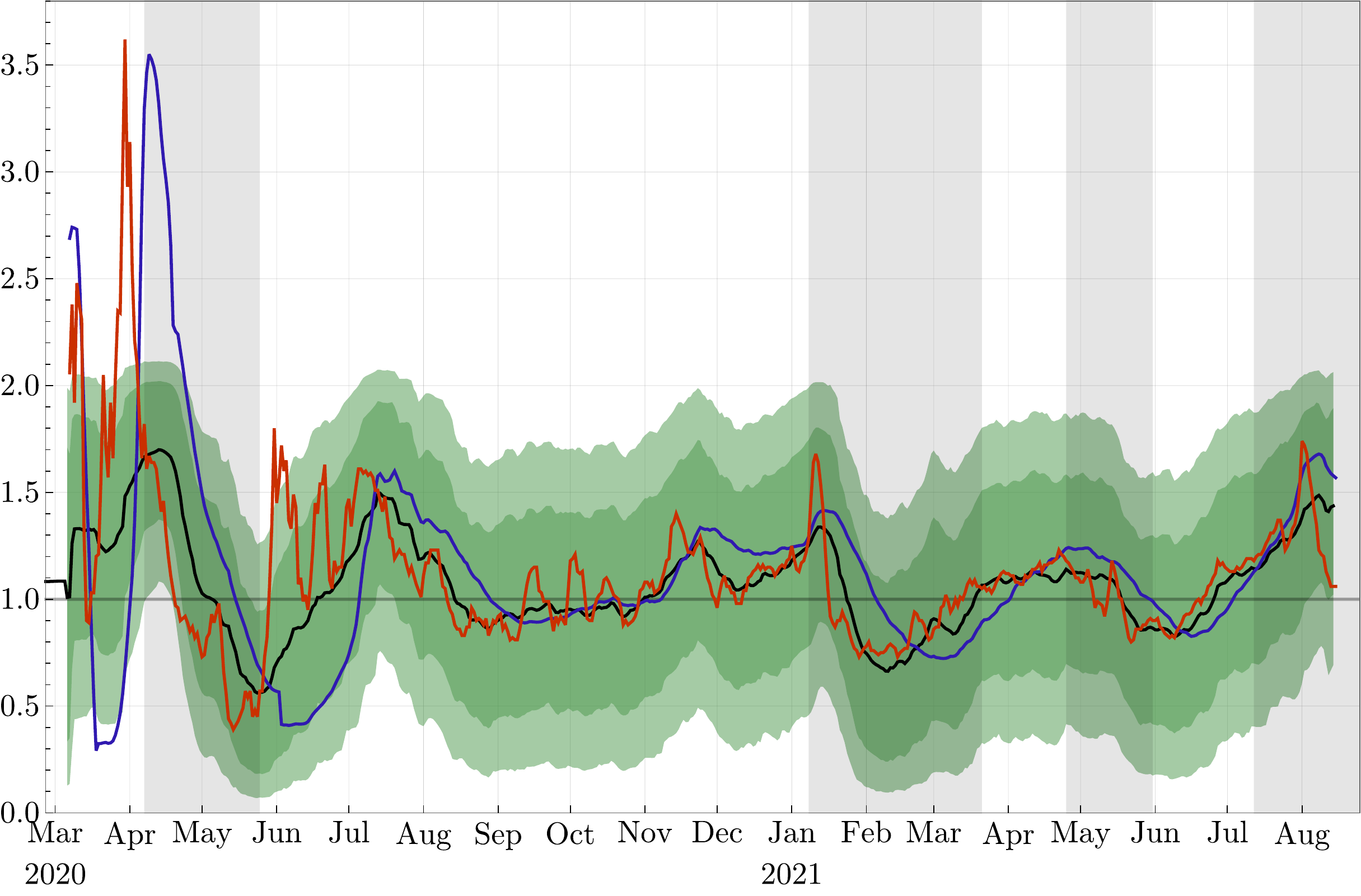}
    \caption{{\it The effective reproduction number for Tokyo computed with different approaches: the black curve and the green confidence intervals are from our approach, the blue curve is provided by \cite{Qiwen}, the red curve is borrowed
                from \cite{Toyo}}}
    \label{fig_comp}
\end{figure}

In Figure \ref{fig_comp} we compare our evaluation for $\rt$ with two other resources. 
One of them is provided by Toyokeizai, a book and magazine publisher
based in Tokyo \cite{Toyo}. Their approach for the computation of $\rt$ is completely different from ours, and is based on a simplified version of a formula proposed
by H.~Nishiura \cite{Nishiura}.
The second resource is a companion paper \cite{Qiwen} in which similar investigations
are performed with a different approach: a continuous model is used for the simulation, 
and the data assimilation part is based on the ensemble Kalman filter. 
Clearly, in the first several months of the epidemic, these different methods lead
to rather diverse values of $\rt$. This phenomena is expected to be due to inaccurate data, irregular release from $H$, and small numbers of agents. On the other hand, since mid-July 2020, the two curves from other studies and our mean $\rt$ share very similar patterns. 
From the figure, one can notice than our approach and the approach of \cite{Qiwen} both provide a more stable $\rt$ estimation than the approach of \cite{Toyo}. Indeed, \cite{Toyo} uses the observations for the computation of $\rt$ directly, while our approach and the approach of \cite{Qiwen} model that pandemic and use data assimilation to produce $\rt$ estimations. As a result, the rapid variations in $\rt$ estimations do not appear anymore.

We believe that our agent-based model, together with the particle filter approach, constitutes a rather simple method for studying the evolution of epidemics.
However, it implicitly takes some assumptions into account, which might not always be satisfied.
For example, we have assumed that the number of agents in $S$ is much
larger than the cumulative number of agents who got infected. This is still true
for Tokyo, with a population of about 14 millions, and a cumulated number of agents
discharged from $H$ reaching a value around $0.3$ million by the end of August 2021.
If the ratio of (recovered agents)\;\!/(susceptible agents) is no longer negligible, it would be necessary to consider the compartment $S$ of a limited size.

Another drastic simplification is the absence of graph structure in our model.
Indeed, we could have considered each agent located on a node of a large graph
(for example homogeneous, random, or temporal)
and implemented the interactions between agents by the edges of the graph.
Such additional structures would have increased the complexity of the model,
but is probably not necessary so far, given the current size of the epidemic in Tokyo.
However, in the future or for a smaller susceptible population (or a bigger infected population), it might be necessary to consider them. Note however that 
though the model would have to be remade, the approach with the particle filter should still be valid.

Finally, new medical inputs could be implemented in our model: so far, the different
variants of the virus have been disregarded, and the ongoing vaccination campaign
is not considered. For an efficient implementation of these new
factors, some of the current compartments should be subdivided or refined.
However, the drawback of such a process would be the increase of the number of parameters, either
provided by some medical services, or which have to be computed by the experiments.
This could unexpectedly enlarge the uncertainty of the predictions made by the experiment, either by the errors of the parameters used, or by the excessive freedoms of the model.
Our future plan is to work on these issues.


\begin{thebibliography}{10}

    \bibitem{Arm}
    E. Armstrong, M. Runge, J. Gerardin,
    \emph{Identifying the measurements required to estimate rates of COVID-19 transmission, infection, and detection, using variational data assimilation},
    Infectious Disease Modelling 6, 133--147, 2021.

    \bibitem{AM}
    F. Arroyo-Marioli, F. Bullano, S. Kucinskas, C. Rond\' on-Moreno,
    \emph{Tracking R of COVID-19: A new real-time
        estimation using the Kalman filter},
    PLoS ONE 16(1): e0244474, 2021.

    \bibitem{BEC}
    D. Buitrago-Garcia, D. Egli-Gany, M. J. Counotte, S. Hossmann, H. Imeri, A.M. Ipekci, et al.,
    \emph{Occurrence and transmission potential of asymptomatic and presymptomatic SARS-CoV-2 infections: A living systematic review and meta-analysis}, PLoS Med 17(9): e1003346, 2020.

    \bibitem{BCB}
    O. Byambasuren, M. Cardona, K. Bell, J. Clark, M.-L. McLaws, P. Glasziou,
    \emph{Estimating the extent of asymptomatic COVID-19 and its potential for community transmission: systematic review and meta-analysis},
    J. Association of Medical Microbiology and Infectious Disease Canada
    5 Issue 4, 223--234, 2020.

    \bibitem{Cheng}
    S. Cheng, C.C. Pain, Y.-K. Guo, R. Arcucci,
    \emph{Real-time updating of dynamic social networks for COVID-19 vaccination strategies}, Preprint \url{https://arxiv.org/abs/2103.00485}.

    \bibitem{Daza}
    M.L. Daza-Torres, M.A. Capistr\' an, A. Capella, J.A. Christen,
    \emph{Bayesian sequential data assimilation for COVID-19 forecasting},
    Preprint \url{https://arxiv.org/abs/2103.06152}.

    \bibitem{EN}
    R. Engbert, M.M. Rabe, R. Kliegl, S. Reich,
    \emph{Sequential Data Assimilation of the Stochastic SEIR
        Epidemic Model for Regional COVID-19 Dynamics},
    Bulletin of Mathematical Biology 83:1, 2021.

    \bibitem{E}
    G. Evensen, J. Amezcua, M. Bocquet, A. Carrassi, A. Farchi, A. Fowler, P.L. Houtekamer,
    C.K. Jones, R.J. de Moraes, M. Pulido, C. Sampson, F.C. Vossepoel,
    \emph{An international initiative of predicting the SARS-CoV-2 pandemic using ensemble data assimilation},  Foundations of Data Science,
    American Institute of Mathematical Sciences, 2020.

    \bibitem{GH}
    R. Ghostine, M. Gharamti, S. Hassrouny, I. Hoteit,
    \emph{ An Extended SEIR Model with Vaccination for
        Forecasting the COVID-19 Pandemic in Saudi Arabia Using an Ensemble
        Kalman Filter}, Mathematics 9, 636, 2021.

    \bibitem{GMc}
    K.M. Gostic, L. McGough , E.B. Baskerville, S. Abbott, K. Joshi, C. Tedijanto, et al.,
    \emph{Practical considerations for measuring the effective
        reproductive number, $R_t$},
    PLoS Comput Biol 16 No. 12,  e1008409, 21 pages, 2020.

    \bibitem{He}
    Z. He, L. Ren, J. Yang, L. Guo, L.  Feng, C. Ma, et al.,
    \emph{Seroprevalence and humoral immune durability of anti-SARS-CoV-2 antibodies in Wuhan, China: a longitudinal, population-level, cross-sectional study},
    The Lancet 397, Issue 10279, 1075--1084, 2021.

    \bibitem{LZL}
    X. Li, Z. Zhao, F. Liu,
    \emph{Big data assimilation to improve the predictability of COVID-19},
    Geography and Sustainability 1(4), 317--320, 2020.

    \bibitem{MI}
    L. Mitchell, A. Arnold,
    \emph{Analyzing the effects of observation function selection in ensemble Kalman filtering for epidemic models},
    Mathematical Biosciences 339, 2021.

    \bibitem{M3}
    m3: \ \url{https://www.m3.com/open/iryoIshin/article/849820/}

    \bibitem{NA}
    P. Nadler,  S. Wang,  R. Arcucci, X. Yang, Y. Guo,
    \emph{An epidemiological modelling approach for COVID-19 via data
        assimilation},
    European Journal of Epidemiology 35, 749--761, 2020.

    \bibitem{Nishiura}
    H. Nishiura: \ \url{https://github.com/contactmodel/COVID19-Japan-Reff}

    \bibitem{OP}
    Osaka prefecture government,
    \emph{Citizens awareness and behavior change of measures against COVID-19},
    \url{http://www.pref.osaka.lg.jp/hodo/attach/hodo-40479\_4.pdf}

    \bibitem{PL}
    A.M. Pollock, J. Lancaster,
    \emph{Asymptomatic transmission of covid-19},
    BMJ, 371:m4851, 2020.

    \bibitem{RE}
    T.C. Rebollo,  D. Coronil,
    \emph{Predictive data assimilation through Reduced Order Modeling for epidemics with data uncertainty},
    Preprint \url{https://arxiv.org/abs/2004.12341}.

    \bibitem{RH}
    C.J. Rhodes, T.D. Hollingsworth,
    \emph{Variational data assimilation with epidemic models},
    Journal of Theoretical Biology 258, 591--602, 2009.

    \bibitem{Sil}
    V.L.S. Silva, C.E. Heaney, Y. Li, C.C. Pain,
    \emph{Data assimilation predictive GAN (DA-PredGAN): applied to determine the spread of COVID-19},
    Preprint  	\url{https://arxiv.org/abs/2105.07729}.

    \bibitem{Qiwen}
    Q. Sun, S. Richard, T. Miyoshi,
    \emph{Analysis of COVID-19 in Japan with Extended SEIR model and ensemble Kalman filter}, in preparation, August 2021.

\bibitem{TMG}
Tokyo Metropolitan Government: \url{https://www.bousai.metro.tokyo.lg.jp/_res/projects/default_project/_page_/001/010/030/2020080608.pdf}

    \bibitem{Toyo}
    Toyokeizai: \ \url{https://toyokeizai.net/sp/visual/tko/covid19/en.html}

    \bibitem{vW}
    J.-D. van Wees, S. Osinga, M. van der Kuip, M. Tanck, M. Hanegraaf,
    M. Pluymaekers, et al.,
    \emph{Forecasting hospitalization and ICU rates of the COVID-19 outbreak: an efficient SEIR model}, Bull World Health Organ. E-pub: 30 March 2020.

    \bibitem{Xu}
    X.-K. Xu, X.F. Liu, Y. Wu, S. Taslim Ali, Z. Du, P. Bosetti, E.H.Y. Lau, B.J. Cowling, L. Wang,
    \emph{Reconstruction of Transmission Pairs for Novel Coronavirus Disease 2019 (COVID-19) in Mainland China: Estimation of Superspreading Events, Serial Interval, and Hazard of Infection},
    Clinical Infectious Diseases 71, Issue 12, 3163--3167, 2020.

\end{thebibliography}
\end{document}